\documentclass[twocolumn]{autart}    

\usepackage{graphicx}
\usepackage[dvipsnames]{xcolor}
\graphicspath{{./fig/}}
\usepackage{amsmath}
\usepackage{amssymb}
\usepackage{mathrsfs}
\usepackage{amsfonts,bbm}
\usepackage{dsfont}
\usepackage{subfigure}
\usepackage{epstopdf}
\usepackage{subfigure}

\usepackage{balance}

\newtheorem{theorem}{Theorem}

\newtheorem{proposition}{Proposition}
\newtheorem{lemma}{Lemma}

\newtheorem{example}{Example}
\newtheorem{remark}{Remark}
\newtheorem{problem}{Problem}
\newtheorem{assumption}{Assumption}
\newtheorem{standing}{Standing Assumption}

\newcommand{\real}{\mathbb{R}}

\newcommand{\R}{\mathbb{R}}

\newcommand{\N}{\mathbb{N}}

\newcommand{\mc}{\mathcal}

\newcommand{\subscr}[2]{{#1}_{\textup{#2}}}
\newcommand{\until}[1]{\{1,\dots,#1\}}

\newcommand*\diag[0]{\mbox{diag}} 

\DeclareSymbolFont{bbold}{U}{bbold}{m}{n}
\DeclareSymbolFontAlphabet{\mathbbold}{bbold}

\newcommand{\vect}[1]{\mathbbold{#1}}

\newcommand{\vectorones}[1][]{\vect{1}_{#1}}
\newcommand{\vectorzeros}[1][]{\vect{0}_{#1}}

\newcommand\oprocendsymbol{\hbox{$\square$}}
\newcommand\oprocend{\relax\ifmmode\else\unskip\hfill\fi\oprocendsymbol}

\begin{document}

\begin{frontmatter}


\title{Gather-and-broadcast frequency control in power systems}

\thanks[footnoteinfo]{A preliminary version of part of the results in this paper is in \cite{FD-SG:15}. 
This material is supported by ETH Z\"urich funds and the SNF Assistant Professor Energy Grant \#160573.}
\author[Paestum]{Florian D\"orfler}\ead{dorfler@control.ee.ethz.ch}    
\quad \quad
\author[Pisa]{Sergio Grammatico}\ead{s.grammatico@tue.nl}               

\address[Paestum]{Automatic Control Laboratory, Swiss Federal Institute of Technology (ETH) Z\"urich, Switzerland, 8092}  
\address[Pisa]{Department of Electrical Engineering, Eindhoven University of Technology, The Netherlands, 5612}             



\begin{abstract}                          
We propose a novel frequency control approach in between centralized and distributed architectures, that is a  continuous-time feedback control version of the dual decomposition optimization method. Specifically, a convex combination of the frequency measurements is centrally aggregated, followed by an integral control and a broadcast signal, which is then optimally allocated at local generation units. 
We show that our {\em gather-and-broadcast} control architecture
comprises many previously proposed strategies as special cases.
We prove local asymptotic stability of the closed-loop equilibria of the considered power system model, which is a nonlinear differential-algebraic system that includes traditional generators, frequency-responsive devices, as well as passive loads, where the sources are already equipped with primary droop control. 
Our feedback control is designed such that the closed-loop equilibria of the power system solve the optimal economic dispatch problem.
\end{abstract}

\end{frontmatter}


\section{Introduction}


The quintessential task of power system operation is to match electrical load and generation. The power balance in an AC power network can be directly accessed via the system frequency, making frequency regulation the fundamental mechanism to ensure the load-generation balance. This task is subject to operational constraints, system stability, and economic interests, and it is traditionally accomplished by adjusting generation in a hierarchical structure consisting of three layers: primary droop control, secondary automatic generation control (AGC), and tertiary control (economic dispatch). These layers range from fast to slow timescales, and from decentralized to centralized control architectures \cite{JM-JWB-JRB:08,AJW-BFW:96}.


With the increasing integration of variable renewable sources, such as wind and solar power, low-inertia power electronic generation, larger peak loads, such as electric vehicles, and liberalized reserve markets on increasingly slower times (and their accompanying deterministic frequency errors), power grids are subject to larger and faster fluctuations \cite{MM-BF-BK-MS:15}. 
In such a distributed generation environment, frequency control requires more fast-ramping generators to act as spinning reserves nowadays mostly provided by gas-driven generation, which is expensive, inefficient, and the resulting emissions defeat the purpose of renewables \cite{LW-KM:02}. As a partial remedy, distributed frequency control through inverter-interfaced sources \cite{JMC-LGF-JTB-EG-RCPG-MMP-JIL-NMA:06} or loads \cite{short:07} has a high potential due to the fast ramping capabilities of these devices. In any case, the task of frequency regulation will have to be shouldered by more and more small-scale and distributed devices. 

From a control perspective, the  main objective of frequency control is to regulate the system frequency subject to operational constraints and economic interests such as load sharing,  optimal generation dispatch, or according to the outcome of reserve markets. Further constraints include a partial information structure accounting for distributed generation, liberalized markets, and limited system knowledge. A plethora of strategies has been developed to address these tasks ranging from fully decentralized to centralized architectures, partially relying on time-scale separation and hierarchical control, and being dependent on the detailed system model, load and generation forecasts. 
%
While centralized strategies such as AGC often suffer from a single point of failure, distributed or fully decentralized approaches 
often fall short in practical implementations and typically require a retrofitting of a {costly peer-to-peer} communication architecture. We postpone a detailed literature review to Section~\ref{subsec: review}, 
where we also present some novel results of independent interest concerning robustness and fairness issues. 

In this paper, we consider a nonlinear, differential-algebraic equation (DAE), and heterogeneous power system model including traditional generation, power electronic sources, and frequency-responsive as well as passive loads. We assume that the sources are already equipped with primary droop control, and we focus on designing the secondary control strategy while simultaneously solving a tertiary economic dispatch problem.
Our control approach falls square in between centralized and distributed architectures, and it is motivated and developed by exploiting parallels in dual decomposition methods in optimization \cite{boyd:admm}, auctions in markets \cite{HRV-JR:10}, mean field control \cite{grammatico:parise:colombino:lygeros:15}, as well as classic AGC \cite{JM-JWB-JRB:08}. 
Interestingly, our 
control architecture
includes many previous frequency control strategies for specific parameter sets. 

Specifically, we start with an online optimization routine for the steady-state dynamics based on {the} dual {decomposition method} that evaluates the price of frequency violation in feedback with the optimal generation response of each generator. Our iterative algorithm resembles a decentralized auction mechanism for a spot market.  Next, we propose a continuous-time feedback control version of this optimization scheme as an aggregation of a convex combination of {frequency} measurements, followed by integral control and optimal local allocations of a broadcast control signal. Our {gather-and-broadcast controller} is such that the closed-loop equilibria of the power system {are optimizers of the} economic dispatch. 
{We believe that our gather-and-broadcast control strategy combines appealing features from both centralized and distributed strategies. It robustifies the frequency control by drawing upon the information of multiple sensors and distributing the control actions to multiple generators,} it does not require any model knowledge, 
{it relies on unidirectional broadcast communication}, and it is privacy preserving: no participant needs to communicate its internal model or cost {function}. 
We prove local asymptotic stability of the nonlinear closed-loop DAE system for a specific class of strictly convex cost functions that give rise to typical secondary control curves encountered in practice, including dead-bands, linear response regions, and saturation effects. The main technical results in this paper
generalize those in our preliminary work \cite{FD-SG:15}, which are based on quadratic cost functions and more restrictive assumptions on the system parameters.
Our analysis relies on a dissipative Hamiltonian formulation of the closed loop {system}, an incremental Bregman-type Lyapunov function as in \cite{ST-MB-CDP:16}, convex analysis \cite{rockafellar:wets}, and a LaSalle invariance principle for DAE systems \cite{FD-JS:ECC16,DJH-IMYM-90}.
%

The paper is organized as follows. In Section \ref{sec:problem} we introduce the frequency control problem that includes both frequency regulation and optimal economic dispatch, and we provide a comprehensive literature review.
In Section \ref{sec:dd_frequency_control} we propose our novel frequency control strategy, and in Section \ref{sec:stability} we show local asymptotic stability of {a desirable subset of the} closed-loop equilibria. In Section \ref{sec:simulations}, 
we illustrate the performance of our strategy with a simulation case study on the IEEE39 New England grid and also compare it to other controllers. Section \ref{sec:conclusion} concludes the paper and raises some open questions.

\subsection*{Notation}
$\R$, $\R_{>0}$, $\R_{\geq 0}$, $\R_{<0}$, $\R_{\leq 0}$ denote the set of real, positive real, non-negative, negative and non-positive real numbers, respectively. $A^\top \in \R^{m \times n}$ denotes the transpose of $A \in \R^{n \times m}$. 
Given some matrices $A_1, \ldots, A_N$, $\text{diag}\left( A_1, \ldots, A_N\right)$ denotes the block diagonal matrix with $A_1, \ldots, A_M$ in block diagonal positions. 
Given some functions or scalars $f_1, \ldots, f_N$, we use the vector notation $\boldsymbol{f} := [ f_1, \ldots, f_N ]^\top$ and matrix notation $\boldsymbol{F} := \textup{diag}\left( f_1, \ldots, f_N\right)$, unless differently specified.
$\vectorones[N]$ ($\vectorzeros[N]$) denotes a vector in $\R^N$ with elements all equal to $1$ ($0$). 
Given a function $f: \R^N \rightarrow \R$, the operator $\boldsymbol{\nabla} f( \boldsymbol{\cdot} ): \R^N \rightarrow \R^N $ denotes the gradient $\left[ \textstyle \frac{\partial f}{ \partial x_1 }( \boldsymbol{x} ), \ldots, \frac{\partial f}{ \partial x_N }( \boldsymbol{x} ) \right]^\top$.
{The sum operator, i.e., $\sum_{i}$ or $\sum_{i,j}$, applies to all terms on its right side as in \cite{rockafellar:wets}.}


\section{Frequency control in power systems} \label{sec:problem}

\subsection{Power system model}
Consider a power system modeled as a graph $G = (\mathcal V,\mathcal E)$ with nodes (or buses) $\mathcal V = \until N$ and edges (or branches) $\mathcal E \subseteq \mathcal{V} \times \mathcal{V}$. With each bus $i \in \mathcal V$, we associate 
a harmonic voltage waveform $V_{i} \cos(\omega^{*}t + \theta_{i})$, where $\omega^{*} = 2\pi \cdot f^{*}$ (and $f^{*}=50 \,\textup{Hz}$ or $f^{*} = 60 \,\textup{Hz}$ is the nominal grid frequency). We consider a lossless high-voltage transmission grid {with topology} induced by the sparse susceptance matrix $\tilde{B} \in \R^{N \times N}$. 
We partition the 
{buses as}
$\mathcal V = \mathcal{G} \cup \mathcal{F} \cup \mathcal{P}$ corresponding to synchronous generators $\mathcal G$, buses with frequency-responsive devices $\mathcal F$ (e.g., frequency-sensitive loads or inverter sources performing droop control), and passive buses $\mathcal P$ (e.g., static loads or inverters performing maximum power-point tracking). 
The associated DAE model reads as \cite{DJH-IMYM-90,JM-JWB-JRB:08} 
\begin{subequations}%
\label{eq:model}%
\begin{align}
	\forall i \in \mathcal{G} &:
	M_{i} \ddot \theta_{i} + D_{i} \dot \theta_{i} = P_{i} + u_{i} - \sum\limits_{j \in \mathcal V} B_{i,j}  \sin(\theta_{i} - \theta_{j})
	\label{eq:model -- 1}
	\\
	\forall i \in \mathcal{F} &:
	D_{i} \dot \theta_{i}  = P_{i} + u_{i} - \sum\limits_{j \in \mathcal V} B_{i,j}  \sin(\theta_{i} - \theta_{j})
	\label{eq:model -- 2}
	\\
	\forall i \in \mathcal{P} &:
	0  = P_{i} + u_{i} - \sum\limits_{j \in \mathcal V} B_{i,j}  \sin(\theta_{i} - \theta_{j})
	\label{eq:model -- 3}
\end{align}
\end{subequations}%
where, for all $i \in \mc{V}$, $P_{i} \in \R$ is a constant power injection or demand (positive for sources and negative for loads), $u_{i} \in \mathcal U_{i} =  [\underline{u}_{i}\,,\,\overline{u}_{i}] \subset \R$ is a controllable injection or demand,  and $B_{i,j} := \tilde{B}_{i,j} V_i V_j$ is the effective susceptance for all $i,j \in \mathcal V$. A generator $i \in \mathcal{G}$ is characterized by its rotational inertia $M_{i}>0$ and primary droop control coefficient $D_{i}>0$. A frequency-responsive device $i \in \mathcal{F}$ is characterized by its frequency-sensitivity $D_{i}>0$ (e.g, the droop coefficient for inverters or actively controlled loads, or the damping of a frequency-dependent load). Passive buses (inverters performing power-point tracking and static loads) have no dynamics. Finally, the absence of integral control at node $i \in \mc V$ is modeled by $\mc U_{i} = \{0\}$.
%

{
\begin{remark}[Unmodeled dynamics]
We do not model reactive power and voltage dynamics, as they do not affect the frequency control problem on the considered time scales -- though all of our forthcoming analyses can be extended under a definiteness assumption on the power flow Jacobian; see \cite{CDP-NM-JS-FD:16} for a related analysis.
{\hfill $\square$}
\end{remark}
}

{Finally, we note that the vector field in \eqref{eq:model} is invariant under a rigid rotation of all angles. Accordingly, all equilibria of the power system model \eqref{eq:model} are sets that are invariant under rigid rotations, and all properties such as uniqueness, optimality, and asymptotic stability of equilibria are to be understood modulo rotational symmetry.}

\subsection{Frequency regulation}
Note that if there is a synchronized solution to \eqref{eq:model} satisfying $\dot \theta_{i} = \omega_{\text{sync}} {\in \mathbb{R}}$ for all $i \in \mathcal V$, then by summing up all  steady-state equations \eqref{eq:model}, the {\em synchronous frequency} (relative to $\omega^{*}$) is obtained from the net power balance:
\begin{equation} \label{eq:Omega}
	\omega_{\text{sync}} := \frac{\sum\nolimits_{i \in \mathcal{V}} P_i + u_{i} }{\sum\nolimits_{i \in \mathcal{G} \cup \mathcal{F}} D_{i}}
	\,.
\end{equation}
If transmission losses are taken into account, there would be another strictly negative term on the right-hand side of \eqref{eq:Omega} depending on the steady-state flow pattern. 

Note that in absence of controllable injections $\{ u_{i} \}_{i \in \mc{V}}$ the synchronous frequency $\subscr{\omega}{sync}$ in \eqref{eq:Omega} is determined by the constant power injections $\{ P_{i} \}_{i \in \mc{V}}$ of possibly slow-ramping generation units, fluctuating renewable sources, and unknown loads. We are interested in regulating the frequency deviation \eqref{eq:Omega} to its nominal (zero) value 
by scheduling the controllable injections $\{ u_{i} \}_{i \in \mathcal{V}}$.

\begin{problem}[Frequency regulation]
\label{prob:regulation}
Schedule the injections $\left\{ u_{i} \in \mathcal U_{i} \right\}_{ i \in \mathcal V }$ to balance load and generation, i.e., so that the frequency deviation $\omega_{\text{sync}}$ in \eqref{eq:Omega} is zero.
{$\hfill \square$}
\end{problem}

\begin{remark}[Multi-area systems]
In an interconnected grid, a second objective aside from load-generation balancing (or equivalently frequency regulation) is to balance the net tie-line interchange power $\subscr{P}{net} \in \real$ (positive for an area generation surplus) over a control area. 
Typically both objectives are unified in a single  {\em area control error} \cite{AJW-BFW:96,JM-JWB-JRB:08} that is added to the frequency error in the integral controllers presented in Section \ref{subsec: review}. Our analysis can be extended to this case by appropriately adding $\subscr{P}{net}$ to the frequency control signals. Henceforth, we restrict ourselves to a single-area grid.
{$\hfill \square$}
\end{remark}

\subsection{Centralized and competitive resource allocation}
\label{subsec: ressource allocation}

A basic {\em feasibility condition} to solve Problem \ref{prob:regulation} is that the total power imbalance can be met by the controllable and {\em constrained} injections $\left\{ u_{i} \in \mathcal U_{i} := [\underline{u}_{i}\,,\,\overline{u}_{i}] \right\}_{i \in \mc{V}}$.

\begin{standing}[Feasibility] \label{ass:feasibility}
$$
-\!\sum\limits_{i \in \mathcal{V}} P_i \in 
\sum_{i \in \mathcal{V}} \mathcal{U}_i = 
\sum\limits_{i \in \mathcal{V}} [\underline{u}_{i}\,,\,\overline{u}_{i}].
\vspace{-0.5cm}$$
\hfill $\square$
\end{standing}

If this feasibility condition is met, then there {are} 
many options to schedule the controllable injections $\{ u_{i} \}_{i \in \mathcal{V}}$ to (asymptotically) regulate $\omega_{\text{sync}}$ in \eqref{eq:Omega} to zero. 

Since we are also interested in solving a resource allocation problem, we associate to every controllable injection a cost function to trade off operating costs, emissions, capacities, and other levels of preference. 

\begin{problem}[Optimal economic dispatch]
\label{prob:optimal-allocation}
\em
Schedule the controllable injections to balance load and generation, while minimizing the aggregate operational cost: 
\begin{align}%
\begin{split}%
	\min_{u \in {\R^N}} & \ \sum\nolimits_{i \in \mathcal{V}} J_{i} (u_{i})
	\\
	\mbox{\textup{s.t.} } & \ \sum\nolimits_{i \in \mathcal{V}} P_i + u_{i} = 0
	\,,
\end{split}%
\label{eq:opf}%
\end{align}
where $J_{i}$ is the cost associated with node $i \in \mathcal{V}$.
{$\hfill \square$}
\end{problem}
Throughout the paper, we consider the following standing assumption.
\begin{standing}[Strict convexity] \label{ass:convexity}
For all $i \in \mathcal{V}$, the {cost} function $J_{i}:\, \mathcal U_{i} \to \R$ is strictly convex and continuously differentiable. 
$\hfill \square$
\end{standing}

Note that we directly incorporate the constraints $\left\{  u_{i} \in \mathcal U_{i}  = [\underline{u}_{i}\,,\,\overline{u}_{i}] \right\}_{i \in \mc{V}}$ in the domain of the cost functions $\{ J_{i} \}_{i \in \mc{V}}$. This can be done also in a smooth way, e.g., using barrier functions.
The economic dispatch problem in \eqref{eq:opf} is typically solved on different time scales and, in a longer planning horizon, it can also include binary unit-commitment constraints and inequality constraints penalizing power flows violating thermal constraints. Here we focus on the reserve scheduling problem, where fast ramping generation and controllable loads are dispatched to meet the \textit{real-time} net demand indicated by the frequency deviation in \eqref{eq:Omega}.

Consider now the Lagrangian function associated with the economic dispatch optimization problem in \eqref{eq:opf}, i.e.,
\begin{equation} \label{eq:Lagrangian}
\mc L( \boldsymbol{u},\lambda) := \sum\nolimits_{i \in \mathcal{V}}  J_{i}(u_{i}) - \lambda \left( u_{i} + P_{i} \right),
\end{equation}
where the scalar $\lambda \in \mathbb R$ is the Lagrange multiplier associated with the constraint $\sum\nolimits_{i \in \mathcal{V}} u_{i} + P_{i} = 0$ in \eqref{eq:opf}.
The necessary KKT optimality conditions \cite{bertsekas:NLP} require that 
\begin{equation}
\textstyle	\frac{\partial \mc L( \boldsymbol{u} ,\lambda)}{\partial \boldsymbol{u}} = \vectorzeros[N]
	\;\implies\; J_{i}^{\prime}(u_{i}^{\star}) = \lambda^{\star} \ \forall i \in \mathcal V\,,
	\label{eq: optimality condition}
\end{equation}
where $J_{i}^{\prime}$ is the derivative of $J_i$.
A basic insight from \eqref{eq: optimality condition} is the {\em economic dispatch criterion} \cite{AJW-BFW:96} stating that all marginal utilities must be identical in the unconstrained case: 
\begin{equation} \label{eq: id-marg-costs}
J_{i}^{\prime}(u_{i}^{\star}) = J_{j}^{\prime}(u_{j}^{\star}) \quad \forall \, i,j \in \mathcal{V}\,.
\end{equation}
So far we took the perspective of centralized social welfare optimization. Motivated by a competitive market perspective (in particular, a spot market), let us now consider the {\em utility maximization} ({that is,}
cost minus benefit minimization) of each market participant $i \in \mc V$:
\begin{equation}
	\min\limits_{u_{i} \in \mathcal U_{i} } \mc L( u_i ,\lambda) = \min_{u_{i}\in \mathcal U_{i} } \left( J_{i}(u_{i}) - \lambda \,  u_{i} \right)
	\,.
	\label{eq: cost-benefit}
\end{equation}
Here, $\lambda$ is the {\em nodal price} which is identical for every participant (in this setup neglecting network congestion). The optimal generation as a function of the price is then obtained by $u_{i}(\lambda^{\star}) =  {J_{i}^{\prime}}^{-1}(\lambda^{\star})$. Accordingly, the constraint in \eqref{eq:opf} can then be formulated as the intersection of the aggregated supply bid and demand curves:
\begin{equation}
0
= 
\sum\nolimits_{i \in \mathcal{V}} P_i + u_{i}
=
\sum\nolimits_{i \in \mathcal{V}} P_i + {J_{i}^{\prime}}^{-1}(\lambda^{\star}).
\label{eq: market clearing}
\end{equation} 
The {\em market clearing price} $\lambda^{*}$ is {then obtained from} \eqref{eq: market clearing}.

Independent of a centralized optimization or a game-theoretic market setup, solving Problem \ref{prob:optimal-allocation} also amounts to asymptotically regulating the frequency, that is, solving Problem \ref{prob:regulation}. Frequency regulation is often referred to as secondary control, whereas offline optimization is referred to as tertiary control. As there are no clear boundaries between these two control objectives, there are many solutions available in the literature to solve the optimal economic power dispatch in \eqref{eq:opf} via online frequency regulation. These solutions range from the classic centralized {\em automatic generation control} (AGC) \cite{AJW-BFW:96,JM-JWB-JRB:08} to distributed optimal frequency {regulation algorithms.} We provide a brief review in the following paragraph together with some lemmas and case studies highlighting common pitfalls, which are of independent interest.

\subsection{Critical review of (de)centralized frequency regulation strategies and their common pitfalls} \label{subsec: review}

\subsubsection*{Decentralized secondary integral control}

To regulate steady-state frequency deviations, one may consider simple {\em decentralized integral controllers}, that is, 
\begin{equation}
k_{i} \, \dot \lambda_{i} = -\omega_{i}, \quad u_{i} = \lambda_{i} \quad \forall i \in \mathcal{K} \subseteq \mc{V}\,,
\label{eq: decentralized integral control}
\end{equation}
where $\mc K \subseteq \mc V$ is the set of sites where integral control is applied with control gain $k_{i}>0$. Such decentralized integral controllers nominally regulate the frequency -- even with global stability guarantees  \cite[Theorems 1 and 2]{CZ-EM-FD:15}, but they also induce additional closed-loop equilibria resulting in undesired injection profiles violating load sharing and economic dispatch objectives \cite[Lemma 4.1]{FD-JWSP-FB:14a}. Moreover, it is well known in power systems \cite{JM-JWB-JRB:08} and in control theory \cite{KJA-TH:06}, that multiple decentralized integral controllers may fail to achieve frequency regulation \cite[Theorem 1]{MA-DVD-HS-KHJ:13} and induce instabilities if the frequency measurements are subject to noise; see also the simulation studies \cite[Section V.B]{MA-DVD-HS-KHJ:13}.
In the following, we make this idea precise and show that decentralized integral control subject to heterogeneous measurement biases {actually leads} to the complete absence of synchronous solutions.

\begin{proposition}[\bf Absence of synchronous solutions] 
\label{Proposition: multiple integrators}
Consider the power system in \eqref{eq:model} under decentralized integral control \eqref{eq: decentralized integral control} subject to biased measurements, i.e.,
\begin{equation}
k_{i} \, \dot \lambda_{i} = -\omega_{i} + \eta_{i}, \quad u_{i} = \lambda_{i} \quad \forall i \in \mc{K} \subseteq \mc{V}\,,
\label{eq: decentralized integral control - bias}
\end{equation}
where
$\eta_{i} \in \real$ is the measurement bias for controller $i \in \mc K$. There exists a synchronous solution $\dot{\boldsymbol{\theta}}$ such that $\lim_{t \rightarrow \infty} \dot{\theta}_i(t) = \subscr{\omega}{sync}$ 
for all $i \in \mc{V}$ only if either \\[.2em]
$\quad (i) \ $ $|\mc K| =1$, i.e., there is only one integral controller; or \\[.2em]
$\quad (ii) \ $ there is $\eta \in \real$ so that $\eta_{i} = \eta$ for all $i \in \mc K$, i.e., all measurement biases are identical.\\[.2em]
In both cases, a synchronous solution $\boldsymbol{\theta}(t)$ satisfies $\lim_{t \rightarrow \infty} \dot{\theta}_i(t) = -\eta$ for all $i \in \mc V$, where $\eta$ is the measurement bias.
{\hfill $\square$}
\end{proposition}

\begin{pf}
Observe that if there exists a synchronous solution $\dot{\theta}(t)$ to the closed loop  \eqref{eq:model}, \eqref{eq: decentralized integral control - bias} such that $\lim_{t \rightarrow \infty}\dot\theta_{i}(t) = \subscr{\omega}{sync}$ for all $i \in \mc V$, then by \eqref{eq: decentralized integral control - bias} this solution must also satisfy
$0 = \lim_{t \rightarrow \infty}\dot\theta_{i}(t) + \eta_i = \subscr{\omega}{sync} + \eta_i$ for all $i \in \mathcal{K}$ . Therefore, $\subscr{\omega}{sync} = -\eta_{i}$ for all $i \in \mc K\subseteq\mc V$ and the statement follows.
{\hfill $\blacksquare$}
\end{pf}

The intuition behind this ``instability'' {mechanism} is that different units -- in spite of being coupled -- aim to stabilize different frequencies and thus fight the physical dynamics. Example \ref{Example: noise induced instability} illustrates this phenomenon.

\begin{figure}[!htb]
	\centering{
	\includegraphics[trim = 4mm 2mm 14mm 0mm, clip=true,
	width=.9\columnwidth]{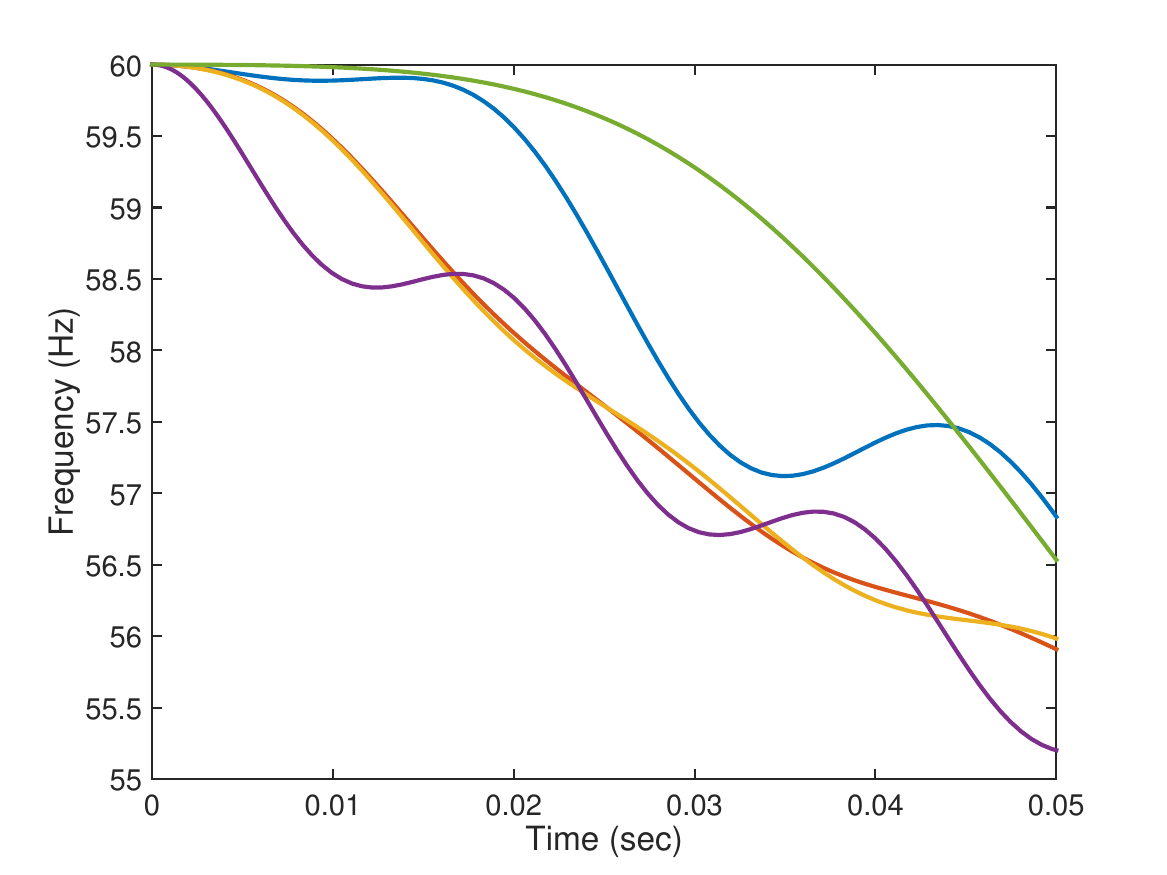}
	\caption{The power system \eqref{eq:model}, with decentralized integral control \eqref{eq: decentralized integral control} is unstable under biased frequency measurements.}
	\label{Fig: noise-induced instability}
	}
\end{figure}

\begin{example}[Bias-induced instability]
\label{Example: noise induced instability}
Consider the power system in \eqref{eq:model} under decentralized integral control \eqref{eq: decentralized integral control - bias}, with measurement biases drawn from a unit variance and zero mean Gaussian distribution. For the system parameters given in Section~\ref{sec:simulations}, the resulting frequency dynamics are unstable as shown in Figure~\ref{Fig: noise-induced instability}. 
\oprocend
\end{example}

\subsubsection*{Automatic generation control}
The industrial standard is the centralized AGC \cite{AJW-BFW:96,JM-JWB-JRB:08} where a frequency measurement is integrated ({together }with the area control error)  at a single site $i^{\star} \in \mathcal{V}$, and the {load-generation} mismatch is allocated to individual generating units according to their {\em participation factors} $\left\{ 1/A_{i}>0 \right\}_{i \in \mathcal{V}}$, often selected as inverse power ratings of the sources, which define their individual contribution:
\begin{equation} \label{eq:AGC 0}
k \, \dot \lambda  = -\omega_{i^{\star}}, \quad u_{i} = \textstyle \frac{1}{A_{i}} \lambda \quad \forall i \in \mathcal{V}.
\end{equation}
Note that the AGC signal in \eqref{eq:AGC 0} may be written as
\begin{equation} \label{eq:AGC}
k \, \dot \lambda  = -\omega_{i^{\star}}, \quad u_{i} = {J_{i}^{\prime}}^{-1}(\lambda) \quad \forall i \in \mathcal{V}\,,
\end{equation}
if the cost function {is} the quadratic function $J_{i}(u_{i}) = \frac{1}{2} A_{i} u_{i}^{2}$. Hence, the AGC strategy \eqref{eq:AGC} achieves identical marginal costs  \eqref{eq: id-marg-costs}
and is implicitly optimal for a quadratic cost.
%
On the other hand, AGC may not be suited for a distributed generation environment due robustness issues (aside from being centralized it relies on a single measurement at $i^{\star}$ that can be compromised) and since a single node $i^{\star}$ may not have the authority to command  the control strategies of all other nodes.

\subsection*{Distributed secondary controllers}
As alternatives to decentralized integral control \eqref{eq: decentralized integral control} or centralized AGC \eqref{eq:AGC 0}, {\em distributed secondary integral controllers} have been proposed that average the integral actions among the generation units through a communication network between the controllers. Different distributed secondary integral approaches have been proposed on the basis of  {\em continuous-time consensus averaging} with all-to-all \cite{HB-JWSP-FD-FB:13j,QS-JG-JMV:13,HL-BJC-WZ-XS:13,MA-DDV-HS-KHJ:14} or nearest-neighbor \cite{JWSP-FD-FB:12u,LYL:13,DB-MDB:14} communication. 
These distributed secondary control approaches can be merged with the tertiary optimization layer, based on the {economic dispatch criterion} \eqref{eq: id-marg-costs} that all marginal utilities must be identical. Different approaches realize this objective based on {continuous-time optimization} approaches  \cite{CZ-UT-SHL:13,FD-JWSP-FB:14a,CZ-EM-FD:15,JZ-FD:14,MA-DVD-HS-KHJ:13,NL-CL-ZC-SHL:13,ST-MB-CDP:16,XZ-AP:14,EM-CZ-SHL:14,SY-LC:14}, game-theoretic ideas \cite{MNE-CAM-NQ:14}, nodal pricing \cite{AJ-ML-VDBPPJ:09}, or discrete-time algorithms \cite{SCT-ADG:12,RM-SD-BBC:12,SK-HG:12}.
All of these algorithms rely on the fact that frequencies should be nominal and marginal costs should be identical in an optimal steady state. 
Such goals are typically achieved by integrating the associated error signals or dual multipliers; see \cite{TS-CDP-AVDS:15} for a 
summary.  
Accordingly, for all $i \in \mathcal{V}$, distributed averaging-based integral (DAI) controllers of the form
\begin{equation}
k_{i} \dot \lambda_{i} = -\omega_{i} + \sum\limits_{j \in \mathcal{V} } w_{i,j} \left( J_{i}^{\prime}(u_{i}) - J_{j}^{\prime}(u_{j}) \right), \;  u_{i} = \lambda_{i},%
\label{eq: DAI}
\end{equation}
are commonly considered, where $W = W^{\top} \in \R^{ N \times N}_{\geq 0}$ induces an undirected and connected communication network. Unlike decentralized integral control \eqref{eq: decentralized integral control}, the distributed averaging-based integral control strategy \eqref{eq: DAI} is robust to measurement errors  \cite[Corollary 4]{MA-DVD-HS-KHJ:13}.

In comparison with the centralized price-based coordination in \eqref{eq: cost-benefit}-\eqref{eq: market clearing} (see also \eqref{eq:AGC}), the strategy in \eqref{eq: DAI} relies on bilateral agreements, and one can imagine scenarios where individual agents aim to maximize their benefit by reporting biased marginal costs to their neighbors (see Example \ref{Example: DAI cheating}).
Another drawback of distributed strategies is that the existing controllable units must be retrofitted with {peer-to-peer} bidirectional communication architecture. 
Finally, aside from the above concerns on operational cost and market power, other issues include the utilities' concern that they give the power system operation out of their hands as well as vulnerabilities to cyber-physical faults and security breaches to which distributed strategies  as in \eqref{eq: DAI} (relying on local sensing and bidirectional communication) may be more susceptible than a broadcast architecture as in AGC \eqref{eq:AGC}.

\begin{figure}[htbp]
	\centering{
	\includegraphics[trim = 18mm 0mm 23mm 0mm, clip=true,width=0.99\columnwidth]{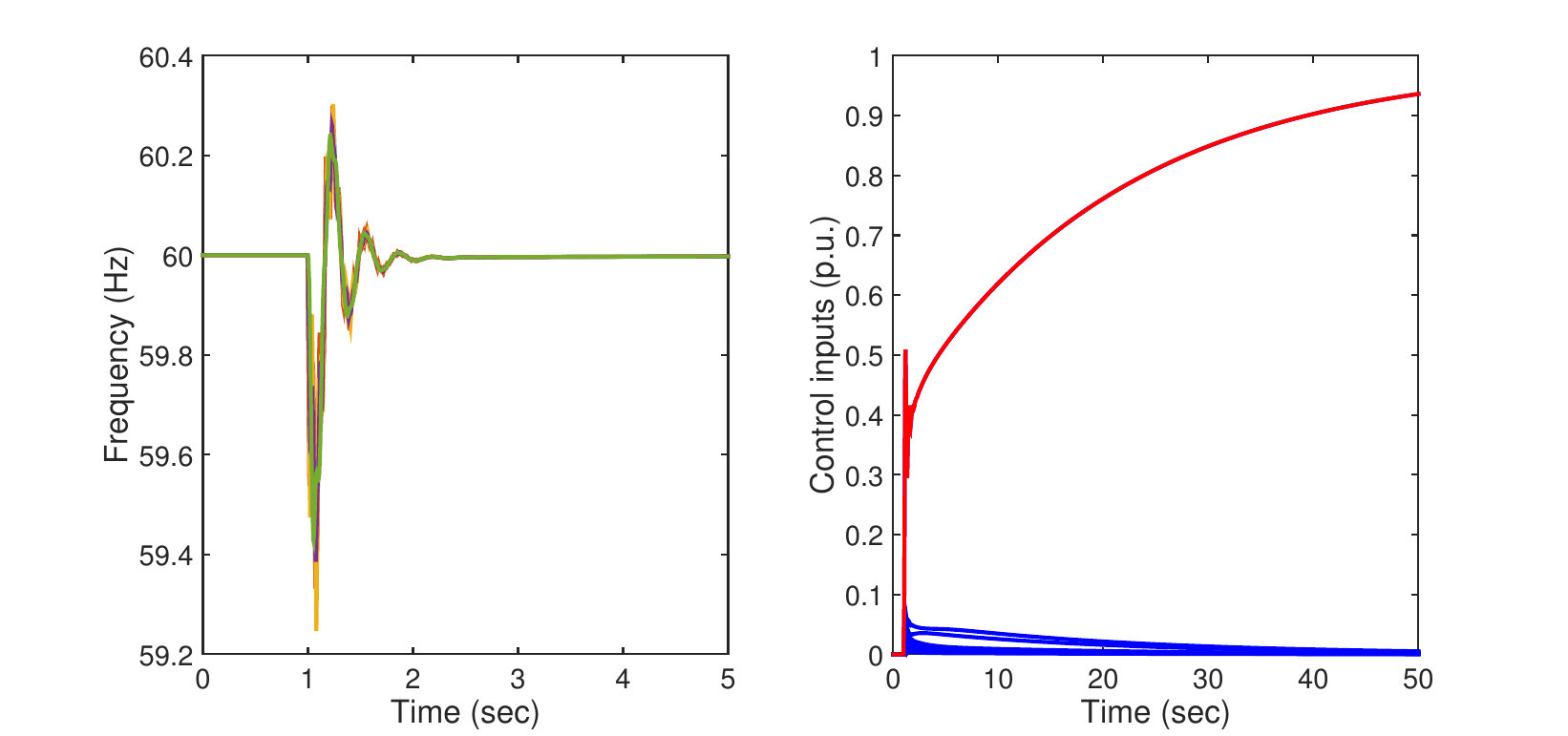}
	\caption{Time-domain plots of frequencies of the closed loop system \eqref{eq:model},\eqref{eq: DAI} and the DAI controller injections if one generator follows the cheating strategy outlined in Example~\ref{Example: DAI cheating}.}
	\label{Fig: cheating}
	}
\end{figure}

\begin{example}[{Cheating under DAI control}]
\label{Example: DAI cheating}
The profit of a unit can be maximized by cheating in the communication protocol of the DAI control \eqref{eq: DAI} as follows. For simplicity, we consider a setup without injection constraints, i.e., $\mc U_{i} = \real$ for all $i \in \mc V$. 
%
Consider the power system \eqref{eq:model} under DAI control \eqref{eq: DAI} and the following strategy executed by node $k$ to maximize its profit:
\begin{enumerate}
	\item[1.] Node $k$ (illegally) reports zero marginal cost to its neighbors so that in all DAI controllers $i \in \mc V \setminus \{k\}$ in \eqref{eq: DAI} we have $ J_{k}^{\prime}(u_{k}) = 0$ as input from node $k$. 

\smallskip

	\item[2.] Node $k$ does not listen to the messages of others nodes by setting $w_{k,j}=0$ for all $j \in \mc V$  in \eqref{eq: DAI}.

\end{enumerate}
%
In this case, the closed-loop system \eqref{eq:model}, \eqref{eq: DAI} admits a steady state with $\dot\theta_{i}=\omega_i=0$ for all $i \in \mc V$, $u_{i}=0$ for all $i \neq k$, and $u_{k} = - \sum_{i \in \mc V}P_{i}$. Hence, node $k$  exclusively regulates the system frequency and accordingly receives a higher compensation than in the fair load sharing case \eqref{eq: id-marg-costs}. Observe, that this modification is not detectable by steady-state frequency measurements.
A simulation of this scenario is shown in Figure~\ref{Fig: cheating} for the system parameters given in Section \ref{sec:simulations}, where generator 10 follows the above cheating strategy to uniquely balance the system and exclusively  receive compensation. We refer to \cite{NM-CDP-JWSP:16} for further cheating mechanisms related to DAI control.
\oprocend
\end{example}

{To conclude this section, let us remark that all integral control mechanisms presented in Subsection \ref{subsec: review} could also have been implemented also as proportional-integral (PI) controllers to improve the system performance and enhance its stability. Indeed, AGC is sometimes implemented through PI control \cite{JM-JWB-JRB:08}. Note that any proportional error feedback is inactive at steady state and does not affect the system equilibria. Thus, even under PI-feedback, the instability mechanism in Proposition \ref{Proposition: multiple integrators} persists, and the cheating strategy in Example \ref{Example: noise induced instability} is still viable.
}



\section{Semi-decentralized frequency control} \label{sec:dd_frequency_control}

\subsection{Market-based discrete-time dual decomposition}

In the following, we provide an alternative frequency control algorithm based on a central price update for violating the net power balance, and inspired by the market-based and game-theoretic insights to the economic dispatch optimization presented in Subsection \ref{subsec: ressource allocation}.

Specifically, we exploit the fact that the Lagrangian function $\mc L(u,\lambda)$ in \eqref{eq:Lagrangian} associated with the optimization problem in \eqref{eq:opf} is separable. Therefore, as the individual costs $\{ J_i \}_{i \in \mathcal{V}}$ are strictly convex and bounded, the optimization problem in \eqref{eq:opf} can be solved iteratively via the {\em gather-and-broadcast} dual decomposition method \cite[Sections 2.1--2.2]{boyd:admm}. This reads as the following iterative {dual-ascent} algorithm, where $k \in \N$ denotes a discrete time index, and $\alpha \in \R_{>0}$ is a sufficiently small step size:
\begin{align}
u_{i}{(k+1)} := & \ \arg\min_{\upsilon \in \mc U_{i}} \left( J_{i}(\upsilon) - \lambda{(k)} \, \upsilon \right) \,, \; \forall i \in \mathcal{V}\,, \label{eq:dual-decomposition-primal} \\
\lambda{(k+1)} := & \ \lambda{(k)} - \alpha \left( \sum\nolimits_{i \in \mathcal{V}} P_i + u_i{(k+1)} \right). \label{eq:dual-decomposition-dual}
\end{align}
At every discrete time step $k$, each node $i \in \mathcal{V}$ computes its optimal injection according to \eqref{eq:dual-decomposition-primal} as a function of the current price $\lambda{(k)}$.
At the same time, the price $\lambda{(k)}$ for the power imbalance is updated in \eqref{eq:dual-decomposition-dual} via a discrete-time integral type control of the power balancing error, that is directly measurable through the frequency signal
\begin{equation} \label{eq:frequency-error}
\omega{({k+1})} := \frac{\sum\nolimits_{i \in \mathcal{V}} P_i + u_i{(k+1)}}{\sum\nolimits_{i \in \mathcal{G} \cup \mathcal{F}} D_{i}} \,.
\end{equation}
Thus, the update \eqref{eq:dual-decomposition-dual} drives the frequency error \eqref{eq:frequency-error} to zero. 
The iteration in \eqref{eq:dual-decomposition-primal}--\eqref{eq:dual-decomposition-dual} can be implemented in a {\em semi-decentralized} fashion. The dual update (or integral control) \eqref{eq:dual-decomposition-dual} determining the current price $\lambda(k)$ is performed at a central site\footnote{In principle, the update in \eqref{eq:dual-decomposition-dual} can also be carried out locally using $\omega_{i}$ instead of \eqref{eq:frequency-error}, since the steady-state frequency error \eqref{eq:frequency-error} is identical throughout the network. However, in a real-time setting such decentralized integral-type updates are subject to the drawbacks listed in Section~\ref{subsec: review}.}%
 using the steady-state frequency error \eqref{eq:frequency-error}, and the primal update \eqref{eq:dual-decomposition-primal} can be carried out locally as a function of the current price $\lambda(k)$ and the cost function $J_{i}(\cdot)$. In this regard, the gather-and-broadcast update in \eqref{eq:dual-decomposition-primal}-\eqref{eq:dual-decomposition-dual} is conceptually similar to AGC \eqref{eq:AGC 0}, with the advantage of guaranteed global convergence \cite[Chapter 6]{bertsekas:NLP}, \cite[Chapter 2]{boyd:admm}, even for non-quadratic costs and local injection constraints.

\begin{proposition}[Discrete-time global convergence]
For a sufficiently small step size $\alpha>0$, the sequence \\ $\textstyle \left( \{u_{i}(k)\}_{i \in \mc{V}}, \, \lambda{(k)} \right)_{k \in \N}$
 defined iteratively in \eqref{eq:dual-decomposition-primal}--\eqref{eq:dual-decomposition-dual} asymptotically converges, from any initial condition, to the unique primal-dual optimal solution to \eqref{eq:opf}. 
$\hfill \square$
\end{proposition}

From a market perspective, the updates \eqref{eq:dual-decomposition-primal} and \eqref{eq:dual-decomposition-dual} correspond to an iterative local utility maximization \eqref{eq:dual-decomposition-primal}, communication of bids $u_{i}^{\star}(k+1)$, subsequent price announcement \eqref{eq:dual-decomposition-dual}, which is again followed by the optimal generation response \eqref{eq:dual-decomposition-primal}, and so on. Such a scheme is referred as an {\em auction} \cite{HRV-JR:10}. Auctions are known to be decentralized yet robust market mechanisms compared to the {\em  exchange trade} based on a (central) price \eqref{eq: cost-benefit}-\eqref{eq: market clearing} and the bilateral {\em over the counter trading} scheme \eqref{eq: DAI} \cite{LH:03} - all of which lead to a Pareto-optimal solution to \eqref{eq: cost-benefit}. 
Finally, since each generation unit follows a best response strategy in \eqref{eq:dual-decomposition-primal}, there are no incentives for unilateral cheating as in Example \ref{Example: DAI cheating} for the DAI controlÊ \eqref{eq: DAI}.

\subsection{Continuous gather-and-broadcast frequency control}
\label{subsec: continuous alg}
In view of the dual decomposition algorithm in \eqref{eq:dual-decomposition-primal}-\eqref{eq:dual-decomposition-dual}, we derive a corresponding {\em continuous-time} version that acts as a feedback control law stabilizing the frequency deviations of the nonlinear DAE model in \eqref{eq:model}. 
We assume that a central aggregator collects a set of frequency measurements in the network and integrates these measurements to form the overall area frequency error as
\begin{equation}\label{eq: mean field}
k \, \dot {\lambda} = -\sum\nolimits_{i \in \mathcal{V}} C_{i} \, \omega_{i},
\end{equation}
where $k>0$ is a gain, and $\{ C_{i}\}_{i \in \mc V} \in [0,1]$ is a set of convex weighting coefficients such that $\sum_{i \in \mc V} C_{i} = 1$. Note that this normalization can be made without loss of generality, through an appropriate scaling of the gain $k$. 
Next the signal $\lambda$ from \eqref{eq: mean field} is broadcast to the individual nodes, where it is dispatched according to 
\begin{equation} \label{eq: allocation}
u_i := {J_i^{\prime}}^{-1}(\lambda) \quad \forall i \in \mc{V}.
\end{equation}
This feedback control scheme relies on the following mean-field-type loop \cite{grammatico:parise:colombino:lygeros:15}: construction of the measurement average $-\sum\nolimits_{i \in \mc V} C_{i} \omega_{i}$ as a global variable that is centrally processed via an integrator and then broadcast back to the individual nodes. 
In the following, we refer to \eqref{eq: mean field}--\eqref{eq: allocation} as {\em gather-and-broadcast} control.
Note that the broadcast-topology is ``one-to-all'', whereas in principle the measurement aggregation can include either only one measurement or possibly all measurements.
Finally, observe that the generation allocation in \eqref{eq: allocation} achieves identical marginal costs $J_{i}^{\prime}(u_{i}^{\star}(t)) = J_{j}^{\prime}(u_{j}^{\star}(t))$ as in \eqref{eq: id-marg-costs} even during transients -- though, feasibility  of \eqref{eq:opf} (i.e., power balance) is  achieved only asymptotically.

\subsection{Comparison with methods proposed in the literature}

For specific parameter choices, the gather-and-broadcast frequency control in \eqref{eq: mean field}--\eqref{eq: allocation} reduces to different control architectures proposed in the literature: 

{\em Automatic generation control} \cite{AJW-BFW:96,JM-JWB-JRB:08}: If only a single measurement coefficient $C_{i}$ is non-zero and each cost function $J_{i}$ is quadratic, then the control scheme in \eqref{eq: mean field}-\eqref{eq: allocation} reduces exactly to the conventional AGC in \eqref{eq:AGC 0}. 

{\em All-to-all  averaging control} \cite{FD-JWSP-FB:14a,HB-JWSP-FD-FB:13j,QS-JG-JMV:13,HL-BJC-WZ-XS:13}:
The gains $C_{i} = D_{i}$ for all $i \in \mathcal{V}$ have been employed for the analysis of centralized averaging-based PI controllers  in  \cite{FD-JWSP-FB:14a,HB-JWSP-FD-FB:13j,MA-DDV-HS-KHJ:14} and experimental implementations in \cite{QS-JG-JMV:13,HL-BJC-WZ-XS:13}.

{\em Mean field control} \cite{grammatico:parise:colombino:lygeros:15}: If all frequencies are weighted equally, $C_i = 1/N$  for all $i \in \mathcal{V}$,  we have a {true mean-field} setup with all nodes treated equally.

{\em Market mechanism} \cite{HRV-JR:10}: The control scheme in \eqref{eq: mean field}--\eqref{eq: allocation} corresponds to an auction mechanism, where the accumulated frequency error 
in \eqref{eq: mean field} serves as pricing signal.

{Note that our controller is inspired by a market mechanism, but in its final form \eqref{eq: mean field}--\eqref{eq: allocation} it does not involve any real-time auctions, prizing, and bidding. Rather it gathers frequency measurements, integrates a convex combination thereof in  \eqref{eq: mean field}, and locally allocates the generation as in \eqref{eq: allocation}. However, extensions are conceivable in the spirit of {\em transactive control} \cite{AA-NT:15} that involve markets and auction mechanisms at discrete and periodic time instants interfaced with the continuous physics through an optimal generation allocation as in \eqref{eq: allocation}.}



\section{Closed-loop stability analysis} \label{sec:stability}

\subsection{Closed-loop equilibria} \label{Subsection: equilibria}
In this section, we analyse the equilibrium states of the closed-loop system in \eqref{eq:model}, \eqref{eq: mean field}--\eqref{eq: allocation}. For a compact presentation, we write the closed-loop system as 
\begin{subequations}%
\label{eq: closed loop}
\begin{align}%
	\dot{\boldsymbol{\theta}}  &= \boldsymbol{\omega}
	\label{eq: closed loop -- 1}\\
	\boldsymbol{M} \dot{\boldsymbol{\omega}}
	&= - \boldsymbol{D} \boldsymbol{\omega} -\boldsymbol{\nabla} U(\boldsymbol{\theta}) + \boldsymbol{P} 
	+ \boldsymbol{{J^{\prime}}^{-1}}(\lambda)
	\label{eq: closed loop -- 2}\\ \smallskip
	k \, {\dot {\lambda}} 
	&= -\boldsymbol{c}^\top \boldsymbol{\omega}
	\label{eq: closed loop -- 3}
	\,,
\end{align}
\end{subequations}
where $\boldsymbol{M} := \diag\left( (M_i)_{i \in \mc{V}}\right)$, $M_i>0$ if $i \in \mc{G}$, $0$ otherwise, 
$\boldsymbol{D} := \diag\left( (D_i)_{i \in \mc{V}}\right)$, $D_i > 0$ if $i \in \mc{G} \cup \mc{F}$, $0$ otherwise, $\boldsymbol{P} := [P_1, \ldots, P_{N}]^\top$; we have introduced the notation
$\boldsymbol{{J^{\prime}}^{-1}}(\cdot) := [ {J_1^{\prime}}^{-1}(\cdot), \ldots, {J_{N}^{\prime}}^{-1}(\cdot) ]^\top: \R \rightarrow \R^{N}$, $\boldsymbol{c} := [ C_1, \ldots, C_{N}]^\top$, and the network potential function 
$U:\, \mathbb{R}^{N} \to \real$ is defined as
\begin{equation}\label{eq: U potential}
\quad U(\boldsymbol{\theta}) := \sum\nolimits_{\{i,j\} \in \mathcal{E}} B_{i,j} \bigl(1- \cos(\theta_{i} - \theta_{j})\bigr).
\end{equation}
Note that $U$ satisfies the overall balance equation $\vectorones[N]^{\top} \boldsymbol{\nabla} U\left( \boldsymbol{\theta}\right) = 0$ due to the symmetry of the power flow.

We remark that when expanding the compact formulation, the equation in \eqref{eq: closed loop -- 2} is still a set of coupled differential and algebraic equations that read as follows:
 \begin{subequations}%
\label{eq: closed loop detail}%
\begin{align}
	\!\forall i \in \mathcal{G} \!&\!:\!&\!
	M_{i} \ddot \theta_{i} + D_{i} \dot \theta_{i} \!&= - \textstyle \frac{\partial U}{\partial \theta_i}( \boldsymbol{\theta} ) + P_{i} +{J_i^{\prime}}^{-1}(\lambda)
	\label{eq: closed loop detail -- 2}
	\\
	\!\forall i \in \mathcal{F} \!&\!:\!&\! 
	D_{i} \dot \theta_{i} \!&= - \textstyle\frac{\partial U}{\partial \theta_i}( \boldsymbol{\theta} ) + P_{i} + {J_i^{\prime}}^{-1}(\lambda)
	\label{eq: closed loop detail -- 3}
	\\
	\!\forall i \in \mathcal{P} \!&\!:\!&\! 
	0 \!&= - \textstyle \frac{\partial U}{\partial \theta_i}( \boldsymbol{\theta} ) + P_{i} + {J_i^{\prime}}^{-1}(\lambda){.}
	\label{eq: closed loop detail -- 4}
\end{align}
\end{subequations}%
The overall closed-loop system in \eqref{eq: closed loop} has the property that if an equilibrium exists, then there exists a unique scalar $\lambda^* \in \R$ such that the power balance is satisfied in steady state. We formalize such a property of the closed-loop equilibria in the following statement.

\begin{proposition}[Closed-loop equilibria]\label{Proposition: equilibria}
The equilibria $\left(\boldsymbol{\theta^*}, \boldsymbol{\omega^*}, \lambda^* \right)$ of the closed-loop system \eqref{eq: closed loop}  are such that $\boldsymbol{\omega^*} = \vectorzeros[N]$, and $\lambda^* \in \real$ is the unique solution to 
\begin{equation} \label{eq: ps}
\sum\nolimits_{i \in \mc{V}} P_i + {J_i^{\prime}}^{-1}(\lambda^*) = 0.
\end{equation}
Moreover, each equilibrium state is an optimal solution to the economic dispatch problem in \eqref{eq:opf}.
$\hfill \square$
\end{proposition}

\begin{pf}
In steady state, we have $\dot{\boldsymbol{\theta}} = \vectorzeros[N]$, $\dot{\boldsymbol{\omega}}=\vectorzeros[N]$, and $\dot{{\lambda}} = 0$. Hence, from \eqref{eq: closed loop -- 1} in steady state we get $\boldsymbol{\omega^*} = \vectorzeros[N] $. Then Equation \eqref{eq: closed loop -- 2} reads in steady state as
\begin{equation}
\label{eq: P steady state}
-\boldsymbol{\nabla} U(\boldsymbol{\theta^{*}}) + \boldsymbol{P} +  \boldsymbol{{J^{\prime}}^{-1}}(\lambda^*)
 = \vectorzeros[N].
\end{equation}
If we multiply this equation from the left by $\vectorones[N]^\top$, since $\vectorones[N]^{\top} \boldsymbol{\nabla} U( \boldsymbol{\theta^{*}})  = 0$, then we obtain Equation \eqref{eq: ps}. 
Since the functions $\{ {J_i^{\prime}}^{-1}(\lambda^*) \}_{i \in \mc{V}}$ are strictly increasing (due to strict convexity), so is the sum $\sum_{i \in \mc{V}} {J_i^{\prime}}^{-1}(\lambda^*)$. Thus, Equation \eqref{eq: ps} admits a unique solution $\lambda^*$. Finally, optimality of the steady-state injections $\{ {J_i^{\prime}}^{-1}(\lambda^*) \}_{i \in \mc{V}}$ follows by construction of the control law in \eqref{eq: allocation}.
{\hfill $\blacksquare$}
\end{pf} 

\subsection{Local asymptotic stability} \label{Subsection: Convergence}
We now perform a stability analysis of the equilibria of the differential-algebraic closed-loop system in \eqref{eq: closed loop}. For simplicity we make use of following assumption.

\begin{assumption}[Scaled cost functions] \label{ass: coeff}
There exists a strictly convex, continuously differentiable function $J: [ \underline{u}, \overline{u}] \rightarrow \R$, for some $\underline{u} \in \R_{<0}$, $\overline{u} \in \R_{>0}$, such that: 
$J^{\prime}(0) = 0$; for all $i \in \mc{V}$, $J_i(\cdot) = J\left( \frac{1}{C_i} \cdot \right)$, where $\{C_i\}_{i \in \mc{V}}$ are as in \eqref{eq: mean field} and positive; $\lim_{u \rightarrow \overline{u}} J(u) = \lim_{u \rightarrow \underline{u}} J(u) = \infty$.
{\hfill $\square$}
\end{assumption}


Assumption \ref{ass: coeff} restricts our frequency control problem to systems with ($i$) each frequency being accounted for in the controller, ($ii$) strictly convex cost functions that are identical up to heterogeneous scaling factors, and ($iii$) frequency measurements that are weighted according to these factors. 
 %
Equivalently, from the perspective of secondary control, it follows from Assumption \ref{ass: coeff} that, for all $i \in \mc{V}$, ${J_i^{\prime}}^{-1}(\cdot) = C_i \, {J^{\prime}}^{-1}(\cdot)$, which implies that 
\begin{equation}\label{eq: c vector}
\boldsymbol{{J^{\prime}}^{-1}}(\cdot) = \left[ \begin{smallmatrix} {J_1^{\prime}}^{-1}(\cdot) \\ \vdots \\ {J_{N}^{\prime}}^{-1}(\cdot) \end{smallmatrix}\right] = \left[ \begin{smallmatrix} C_1 \\ \vdots \\ C_{N} \end{smallmatrix}\right] {J^{\prime}}^{-1}(\cdot) = \boldsymbol{c} \, {J^{\prime}}^{-1}(\cdot).
\end{equation}
Consequently, each unit in \eqref{eq: closed loop} applies the same control input up to a scaling factor. This includes the usual linear controllers but also more general strategies such as commonly encountered frequency response curves with linear regions, deadband (around the nominal frequency), and saturation (at the capacity of the  unit) \cite{JM-JWB-JRB:08}. These curves are scaled for each unit similarly as primary control curves are scaled according to the bid capacity. 

We now state the main technical result of the paper, that is, the local asymptotic stability of the equilibria of the differential-algebraic nonlinear closed-loop system  \eqref{eq: closed loop}, {asymptotic} frequency regulation (Problem \ref{prob:regulation}) {and} the optimal economic power dispatch (Problem \ref{prob:optimal-allocation}).

\begin{theorem}
\label{Theorem: stability}
\textbf{\textup{(Local asymptotic stability and steady-state optimality)}}
If Assumption \ref{ass: coeff} holds, then any equilibrium of the closed-loop system in \eqref{eq: closed loop} satisfying $|\theta_{i}^{*} - \theta_{j}^{*}| < \pi/2$ for all $\{i,j\} \in \mathcal E$ is locally asymptotically stable.
The control inputs $\{ u_{i}(\cdot) \}_{i \in \mc{V}}$ defined in \eqref{eq: allocation} satisfy the optimal dispatch criterion in \eqref{eq: id-marg-costs} for all $t\geq 0$, and asymptotically solve Problems \ref{prob:regulation}, \ref{prob:optimal-allocation}.
{\hfill $\square$}
\end{theorem}

{
The proof of Theorem \ref{Theorem: stability} relies on the construction of an Hamiltonian function which includes the Lur\'e-type integral function $\mathcal{I}: \R \rightarrow \R$ 
\begin{equation} \label{eq: Lure}
	\mathcal{I}(\lambda) := \int_{\lambda_0}^{\lambda} {J^{\prime}}^{-1} \left(\xi \right) \textup{d} \xi \,,
\end{equation}
for some $\lambda_0 \in \R$, characterized by the following properties.
\begin{lemma}
\label{lem: compact I}
The function 
$$ \lambda \mapsto \mathcal{I}(\lambda) - \mathcal{I}(\lambda^*) - \mathcal{I}^{\prime}\left( \lambda - \lambda^*\right) $$
from \eqref{eq: Hamiltonian}, \eqref{eq: Lure} is strictly convex with unique minimizer 
$\lambda = \lambda^*$, radially unbounded and has compact sublevel sets.
{\hfill $\square$}
\end{lemma}

\begin{pf}
First note that $\mathcal{I}$ is continuous, hence locally bounded and with closed sublevel sets 
$\text{lev}_{  \mathcal{I}}\left( \cdot \right) :=  \{ x \in \mathbb{R} \mid \mathcal{I}(x) \leq \cdot \}$ \cite[Theorem 1.6]{rockafellar:wets}.
Then we note that $J$ being strictly convex is equivalent to $J^{\prime}(\cdot)$ being strictly increasing 
\cite[Theorem 12.17]{rockafellar:wets}.
Moreover, it follows from the inverse function theorem \cite[Theorem 2.11]{spivak} that $\left( J^{\prime} \right)^{-1}(\cdot) = \mathcal{I}^{\prime}(\cdot)$ is strictly increasing. In turn, it follows that $\mathcal{I}$ is strictly convex \cite[Theorem 12.17]{rockafellar:wets}.

Since $J$ is continuously differentiable, $\displaystyle \lim_{u \rightarrow \overline{u}} J(u) = \lim_{u \rightarrow \underline{u}} J(u) = \infty$ implies that 
$\displaystyle \lim_{u \rightarrow \overline{u}} J^{\prime}(u) = \infty$ and $\displaystyle \lim_{u \rightarrow \underline{u}} J^{\prime}(u) = -\infty$, and in turn 
$\displaystyle \lim_{y \rightarrow \infty} \left( J^{\prime}\right)^{-1} (y) = \overline{u} > 0$ and $\lim_{y \rightarrow -\infty} \left( J^{\prime}\right)^{-1} (y) = \underline{u} < 0$.
Therefore, $\lim_{\lambda \rightarrow \infty } \mathcal{I}({\lambda}) = \int_{0}^{\infty} \left( J^{\prime} \right)^{-1}( \xi ) \, \textup{d}\xi = \infty$. 
Symmetrically, $\lim_{\lambda \rightarrow -\infty} \mathcal{I}({\lambda}) = -\infty$. This implies that $\mathcal{I}$ is radially unbounded, and thus $\mathcal{I}$ has compact sublevel sets.

Therefore, also the function $\lambda \mapsto \mathcal{I}(\lambda) - \mathcal{I}(\lambda^*) - \mathcal{I}^{\prime}(\lambda^*) \left(\lambda -\lambda^*\right)$ is strictly convex, with a unique minimizer, radially unbounded, and compact sublevel sets. 

We characterize its minimizer via the first order condition, that is, $\mathcal{I}^{\prime}(\lambda) - \mathcal{I}^{\prime}(\lambda^*) = 0$, which is equivalent to ${J^{\prime}}^{-1}(\lambda) = {J^{\prime}}^{-1}(\lambda^*)$. Since the solution must be unique, we conclude that $\lambda = \lambda^*$ is the unique minimizer.
{\hfill $\blacksquare$}
\end{pf}
}

\begin{pf}\textbf{(Theorem \ref{Theorem: stability})}
By employing the steady-state formulation in \eqref{eq: P steady state} and due to the special structure \eqref{eq: c vector}, we can write the closed-loop system in \eqref{eq: closed loop} compactly as follows.
\begin{align}%
	\dot{\boldsymbol{\theta}}  &= \boldsymbol{\omega}
	\nonumber\\
	\boldsymbol{M} \dot{\boldsymbol{\omega}}
	&= - \boldsymbol{D} \boldsymbol{\omega} - \bigl(\boldsymbol{\nabla} U(\boldsymbol{\theta}) - \boldsymbol{\nabla} U(\boldsymbol{\theta^*}) \bigr)  \nonumber \\
    & \quad \, + \left( \boldsymbol{{J^{\prime}}^{-1}}(\lambda) - \boldsymbol{{J^{\prime}}^{-1}}(\lambda^*) \right)
	\nonumber\\ \smallskip
	k \, {\dot {\lambda}} 
	&= -\boldsymbol{c}^\top \boldsymbol{\omega}
	\label{eq: very compact}%
\end{align}
The compact form in \eqref{eq: very compact} reveals the dissipative nature of the closed loop when using the Hamiltonian function
\begin{multline} \label{eq: Hamiltonian}
\mathcal{H}\left( \boldsymbol{\theta}, \boldsymbol{\omega}, \lambda\right) := 
U(\boldsymbol{\theta}) - U( \boldsymbol{\theta^{*}}) - \boldsymbol{\nabla} U( \boldsymbol{\theta^{*}}) \left( \boldsymbol{\theta} - \boldsymbol{\theta^*} \right) + \\
\frac{1}{2} \boldsymbol{\omega}^\top \boldsymbol{M} \boldsymbol{\omega} + {k }\left(
\mathcal{I}(\lambda) - \mathcal{I}(\lambda^{*}) - \mathcal{I}^{\prime}(\lambda^{*}) \left(\lambda-\lambda^{*}\right) \right)
\,,
\end{multline}
{where $\mathcal{I}$ is the Lur\'e-type integral function in \eqref{eq: Lure}.}
In \eqref{eq: Hamiltonian} we used the Bregman distance of $U(\theta)$ and   $\mathcal{I}(\lambda)$ to $\theta^{*}$ and $\lambda^{*}$, respectively, to construct an incremental Hamiltonian function as in \cite{ST-MB-CDP:16}. 
Next we show some properties related to this integral function.

%
To proceed, we calculate the derivative of $\mathcal{H}\left( \boldsymbol{\theta}, \boldsymbol{\omega}, {\lambda}\right)$ along trajectories of the closed-loop system in \eqref{eq: very compact}  as
\begin{equation} \label{eq: Hamiltonian derivative}
\begin{array}{l}
\dot{\mc{H}}\left( \boldsymbol{\theta}, \boldsymbol{\omega}, {\lambda} \right) 
\\
= \left[ 
\begin{matrix}
\boldsymbol{\nabla} U(\boldsymbol{\theta}) - \boldsymbol{\nabla} U(\boldsymbol{\theta^*}) \\ 
\boldsymbol{M \omega} \\ 
k \, \left( \mc{I}^{\prime}(\lambda) - \mc{I}^{\prime}(\lambda^*) \right)
\end{matrix} \right]^\top 
\left[ \begin{matrix} \dot{\boldsymbol{\theta}} \\ \dot{\boldsymbol{\omega}} \\ \dot{\lambda} \end{matrix} \right] \smallskip \\
= \left( \boldsymbol{\nabla} U( \boldsymbol{\theta} ) - \boldsymbol{\nabla} U(\boldsymbol{\theta^*}) \right)^\top \boldsymbol{\omega} + \boldsymbol{\omega}^\top \boldsymbol{M} \dot{\boldsymbol{\omega}} \, + \\
\qquad \left( \mc{I}^{\prime}(\lambda) - \mc{I}^{\prime}(\lambda^*) \right) \, k \, \dot\lambda \\
= \left( \boldsymbol{\nabla} U(\boldsymbol{\theta}) - \boldsymbol{\nabla} U( \boldsymbol{\theta^{*}}) \right)^\top \boldsymbol{\omega} \, + \boldsymbol{\omega}^\top \cdot \\
\ \cdot \left( - \boldsymbol{D} \boldsymbol{\omega} -\! \left(\boldsymbol{\nabla} U(\boldsymbol{\theta}) - \boldsymbol{\nabla} U(\boldsymbol{\theta^{*}}) \right) + \boldsymbol{c} \, \mc{I}^{\prime}(\lambda)  - \boldsymbol{c} \, \mc{I}^{\prime}(\lambda^*) \right) \\
\ \ - \left( \mc{I}^{\prime}(\lambda) - \mc{I}^{\prime}(\lambda^*) \right) \, \boldsymbol{c}^\top \boldsymbol{\omega} \\
= -\boldsymbol{\omega}^\top \boldsymbol{D} \boldsymbol{\omega} \leq 0\,.
\end{array}
\end{equation}
Thus, the Hamiltonian function $\mc{H}$ is non-increasing along the closed-loop trajectories, and we are in a position to apply the LaSalle invariance principle for DAE systems  \cite[Theorem 2.5]{FD-JS:ECC16}, \cite[Theorem 3]{DJH-IMYM-90}.
Namely, we need to construct a compact set ($i$) in which the vector field in \eqref{eq: very compact} is twice continuously differentiable, ($ii$) within which the Hamiltonian has a strict minimum (modulo rotational symmetry) at the desired equilibrium
, and ($iii$) in which the Jacobian of the algebraic equations in \eqref{eq: closed loop detail -- 4} with respect to the algebraic variable $\boldsymbol{\theta}_{\mc P} := [ (\theta_i)_{i \in \mc{P}} ]^\top$ is nonsingular. 
 
In the following, we show that the sublevel set
\begin{multline} \label{eq:Omega_c}
	\Omega_{\rho} := \left\{ \left. (\boldsymbol{\omega},\boldsymbol{\theta},\lambda) \in \real^{2 N +1} \right| \mc{H}( \boldsymbol{\omega}, \boldsymbol{\theta},\lambda) \leq \rho \,, \right.
	\\
	\left. \ |\theta_{i} - \theta_{j}| < \pi/2 \;\;\forall \{i,j\} \in \mathcal E \right\}
\end{multline}
satisfies all of the above conditions, for sufficiently small $\rho>0$.
First, observe that the vector field in \eqref{eq: very compact} is twice continuously differentiable in $\Omega_{\rho}$. 
Next, we show that the dynamics in \eqref{eq: very compact} are bounded in $\Omega_{\rho}$. Note that the Hessian of $U(\cdot)$ has $\{i,j\}$ element $\left[ \boldsymbol{\nabla^2} U(\boldsymbol{\theta}) \right]_{i,j}$ equal to%
\begin{equation}
\label{eq: Hessian}
\frac{ \partial^2 U(\boldsymbol{\theta})}{ \partial \theta_j \partial \theta_i } =
\left\{
\begin{array}{ll}
- B_{i,j} \cos\left( \theta_i - \theta_j \right) & \text{if } j \neq i \\
\sum_{k=1, k\neq i}^{N} B_{i,k} \cos\left( \theta_i - \theta_k \right) & \text{if } j=i.
\end{array}
\right.\!
\end{equation}
%
%
%
By assumption, $|\theta_{i}^{*} - \theta_{j}^{*}| < \pi/2$ for all $\{i,j\} \in \mathcal E$.
Therefore, $\boldsymbol{\nabla^2} U(\boldsymbol{\theta^*})$ is a positive semidefinite and irreducible (due to connectivity) Laplacian matrix {with nullspace aligned along $\vectorones[N]$}
 corresponding to the rotational symmetry, i.e., the dynamics in \eqref{eq: very compact} are invariant under a rigid rotation of all angles $\boldsymbol{\theta}$. Hence,  within $\Omega_{\rho}$ in \eqref{eq:Omega_c}, $U(\boldsymbol{\theta})$ is locally positive definite and its sublevel sets are compact (modulo rotational symmetry). 
Thanks to this fact and due to Lemma~\ref{lem: compact I}, the Hamiltonian function $\mc{H}$ is locally positive definite with respect to an equilibrium $(\boldsymbol{\theta^{*}},\vectorzeros[N],\lambda^*) \in \Omega_{\rho}$, and its sublevel sets are compact (modulo rotational symmetry). 
 
The above reasoning together with the fact that $\dot{\mc{H}}\leq 0$ guarantees boundedness of the frequencies $\boldsymbol{\omega}$, the integral variable $\lambda$, as well as the relative angles $\theta_{i} - \theta_{j}$ for all $\{i,j\} \in \mc{E}$, that is,  
$\boldsymbol{v}^{\top} \boldsymbol{\theta}$ is bounded for any $\boldsymbol{v} \in \R^{N}$ such that $\boldsymbol{v} \perp \vectorones[N]$. To show boundedness of the remaining coordinate, the sum of all angles $\vectorones[N]^{\top} \boldsymbol{\theta}$, we  integrate the controller equation in \eqref{eq: very compact} as $k (\lambda(t) - \lambda(0)) = \boldsymbol{c}^{\top}(\boldsymbol{\theta}(t)-\boldsymbol{\theta}(0))$. Since $\lambda(t)$ is bounded, it follows that $\boldsymbol{c}^{\top} \boldsymbol{\theta}(t) = \vectorones[N]^{\top} \boldsymbol{C} \boldsymbol{\theta}(t)$ is bounded. Thus, both $\vectorones[N]^{\top} \boldsymbol{C} \boldsymbol{\theta}$ and $\boldsymbol{v}^{\top} \boldsymbol{\theta}$ are bounded for any $\boldsymbol{v} \perp \vectorones[N]$. It follows that $\vectorones[N]^{\top} \boldsymbol{\theta}$ is bounded, and the overall dynamics \eqref{eq: closed loop} are bounded. 
 
Finally, within $\Omega_{\rho}$, the Jacobian matrix associated with the algebraic equation \eqref{eq: closed loop detail -- 4}  is a principal submatrix of the irreducible Laplacian matrix in \eqref{eq: Hessian}. Since submatrices of irreducible Laplacians are nonsingular \cite[Lemma 2.1]{FD-FB:11d}, it follows from the implicit function theorem \cite[Theorem 9.28]{rudin-math-analysis} that the algebraic equations \eqref{eq: closed loop detail --  4} are solvable with respect to the algebraic variable $\boldsymbol{\theta}_{\mc P}$.

Therefore, all the conditions of the LaSalle invariance principle for DAE systems are met. It follows that 
the closed-loop trajectories asymptotically converge to largest invariant set in $\Omega_{\rho}$ satisfying $\dot{\mc{H}}(\boldsymbol{\omega},\boldsymbol{\theta},\lambda) = 0$. 


It remains to be shown that these trajectories converge to a desired equilibrium $(\boldsymbol{\theta^{*}},\vectorzeros[N],\lambda^*)$.
From \eqref{eq: Hamiltonian derivative} we conclude that $\dot{\mc{H}}(\boldsymbol{\omega},\boldsymbol{\theta}, \lambda) = 0$ implies $\boldsymbol{\omega}_{\mc{G} \cup \mc{F}} = \vectorzeros[{ |\mc{G}|+|\mc{F}|}]$. The latter implies for the dynamics in \eqref{eq: very compact} that 
\begin{align}
& \bigl( \boldsymbol{\nabla} U(\boldsymbol{\theta}) - \boldsymbol{\nabla} U(\boldsymbol{\theta^{*}}) \bigr) =  \boldsymbol{c} \, \left( {{J^{\prime}}^{-1}}(\lambda) - {{J^{\prime}}^{-1}}(\lambda^*) \right)\,,
\nonumber\\
& k \, \dot \lambda = -\sum\nolimits_{i \in \mc P} C_{i} \, \omega_{i}
\,,
\label{eq: eq cond}
\end{align}
that is, equations \eqref{eq: closed loop detail -- 2}, \eqref{eq: closed loop detail -- 3} are in equilibrium, the algebraic constraint \eqref{eq: closed loop detail -- 4} is met, {and the integral state $\lambda$ is possibly non-stationary so that the control input {$\boldsymbol{u} = \boldsymbol{{J^{\prime}}^{-1}}(\lambda)$} is possibly not at equilibrium $\boldsymbol{{J^{\prime}}^{-1}}(\lambda^{*})$.} 

By summing over the first set of equations in \eqref{eq: eq cond}, we find that $0 = \vectorones[N]^{\top} \, \boldsymbol{c} \cdot \left( {{J^{\prime}}^{-1}}(\lambda) - {{J^{\prime}}^{-1}}(\lambda^*) \right) = {{J^{\prime}}^{-1}}(\lambda) - {{J^{\prime}}^{-1}}(\lambda^*)$, that is, the integral variable $\lambda$ must be at the equilibrium $\lambda^*$. {Thus,} $\bigl( \boldsymbol{\nabla} U(\boldsymbol{\theta}) - \boldsymbol{\nabla} U(\boldsymbol{\theta^{*}}) \bigr) = \vectorzeros[N]$, so that {the system \eqref{eq: very compact} is} at equilibrium (modulo symmetry). 
{ 
Note that the system \eqref{eq: very compact} being at an equilibrium implies that all dynamic states $(\boldsymbol{\theta}_{\mc{G} \cup \mc{F}},\boldsymbol\omega_{\mc{G} \cup \mc{F}},\lambda)$ are stationary, and so all algebraic states $\boldsymbol{\theta}_{\mc P}$  must be stationary as well due to non-singularity of the algebraic equation.
}

It follows that any equilibrium of the closed loop \eqref{eq: closed loop} satisfying $|\theta_{i}^{*} - \theta_{j}^{*}| < \pi/2$ for all $\{i,j\} \in \mathcal E$ is locally asymptotically stable.
The condition in \eqref{eq: id-marg-costs} on the identical marginal costs follows directly by \eqref{eq: allocation}, 
as $J_i\left( u_i(t) \right) = \lambda$ for all $i \in \mc{V}$ and all $t \geq 0$,
which implies that the control inputs $\lim_{t \rightarrow \infty} \{ u_i(t) \}_{i \in \mc V}$ solve Problem \ref{prob:optimal-allocation}.
{\hfill $\blacksquare$}
\end{pf}

Note that the additional assumption $|\theta_{i}^{*} - \theta_{j}^{*}| < \pi/2$ in Theorem \ref{Theorem: stability} is common in power system analysis where it is also referred to as a security or thermal limit constraint and restricts the solution space to desirable power flows.


\section{Numerical simulations} \label{sec:simulations}
\begin{figure}[h!]
\centering
\includegraphics[width = 0.9\columnwidth]{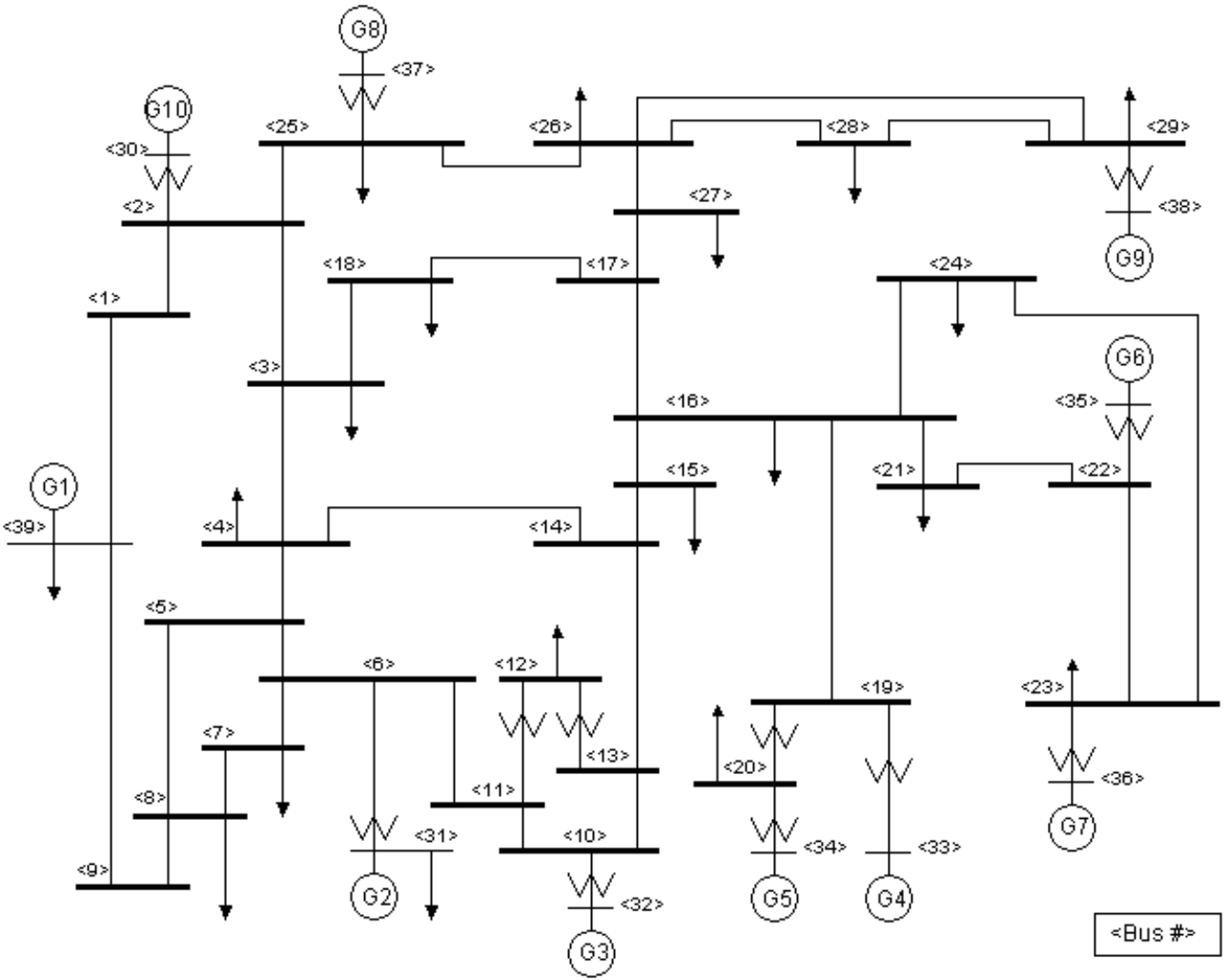}
\caption{IEEE New England test power system.}
\label{fig:IEEE_NE}
\end{figure}
We evaluate the performance of different controllers on the IEEE New England power grid shown in Figure \ref{fig:IEEE_NE}. 
The system has $10$ generators and $39$ buses, serving a total load of about $6$ GW. The generator inertia coefficients $\{M_{i}\}_{i}$ and line susceptances $\{ B_{i,j}\}_{\{i,j\}}$ are obtained from the Power System Toolbox {\cite{JHC-KWC:92}}. The droop coefficients are chosen uniformly $D_i = 1$ for all buses $i$. The cost coefficients $\{ C_i\}_i$ are randomly generated, uniformly in $(0,1)$. We set the integral {gains of all controllers} to $k = 60$.
We simulate {the same scenario as in \cite{CZ-EM-FD:15}:} at time $t=1 \, $s, the demand changes by $33$ MW at buses $4$, $12$ and $20$, creating a power imbalance and causing {the frequencies to drop below the nominal} $60$ Hz.

\begin{figure*}[htbp]
	\centering
	\includegraphics[width=1\textwidth]{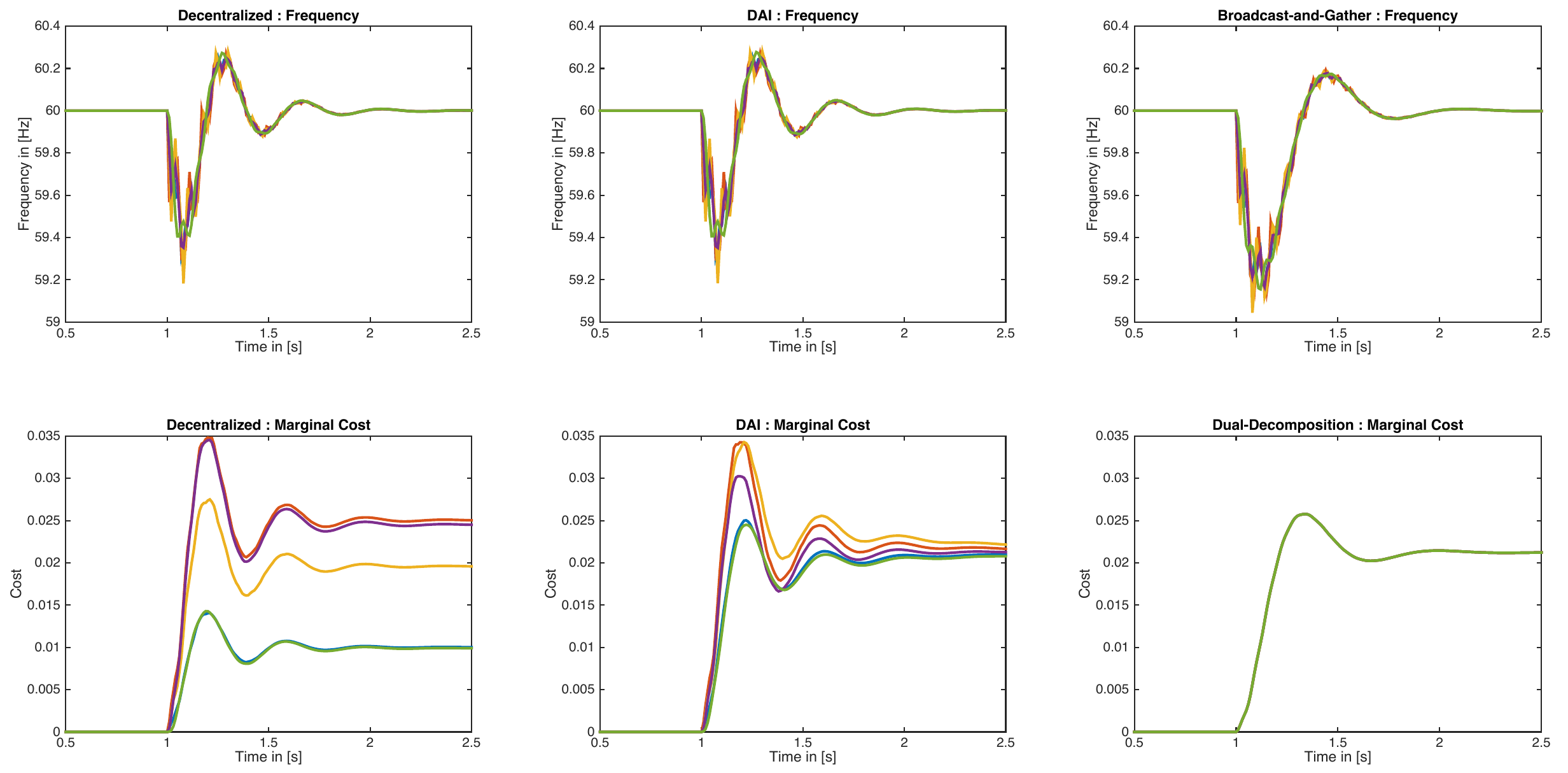}
	\caption{Frequency and marginal costs of generators $2$, $4$, $8$, $10$, for the different control schemes.}
	\label{fig:comparison}
\end{figure*}

We compare our gather-and-broadcast control strategy \eqref{eq: mean field}--\eqref{eq: allocation} for quadratic cost {functions} $J_{i}(\lambda) = \lambda^{2}/C_{i}$ and resulting linear controllers {$u_{i} = C_{i}\, \lambda$} with the fully decentralized {linear} integral control \cite[Section III]{CZ-EM-FD:15} and with the {linear DAI} control \eqref{eq: DAI} based on the same communication graph as in \cite[Section V]{CZ-EM-FD:15}. Figure~\ref{fig:comparison} shows the frequencies and the marginal costs of five generators for the three control schemes. 
Observe that all controllers drive the system frequencies to the nominal value. As expected, the decentralized integral controller does not achieve asymptotically identical marginal costs; on the other hand, both the DAI controller and the proposed one asymptotically solve the optimal economic dispatch problem. In addition, our gather-and-broadcast control guarantees identical marginal costs even during transients. 
{From the qualitative point of view,} Figure \ref{fig:comparison} shows that the closed-loop frequency response induced by our controller is {comparable} to that of the DAI controller, despite {many fewer tuning gains and} the communication requirements being {lower}.

\begin{figure*}[htbp]
	\centering{
	\includegraphics[width=1\textwidth]{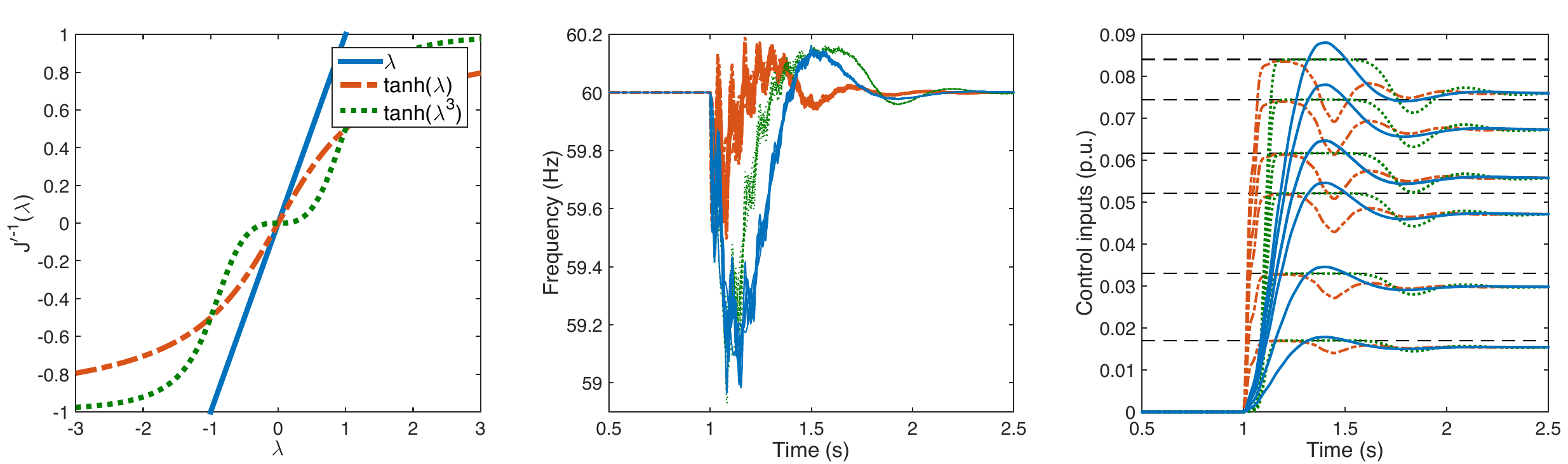}
	\caption{Frequency and control inputs under different nonlinear {controllers} \eqref{eq: nonlinear control} and a linear {one.}}
	\label{Fig: comparison of nonlinearities}
	}
\end{figure*}

In the following, we illustrate the degrees of freedom in the choice of {the} cost functions in the gather-and-broadcast control. Consider the class of controllers {$u_{i} = C_i \, {J^{\prime}}^{-1}(\lambda)$, for $i \in \mc{V}$,}  according to the family {of functions}
\begin{equation}
	{J^{\prime}}^{-1}(\lambda) = \tanh\left(k_{1} \cdot \lambda^{k_{2}}\right)
	\label{eq: nonlinear control}
\end{equation}
{that are} smooth and strictly increasing nonlinear functions encoding a saturation, a smoothly approximated deadzone, and a linear behavior parameterized by $k_{1}>0$ {and an odd $k_{2}>0$.} The associated cost function $J(\cdot)$ satisfies Assumption \ref{ass: coeff}. Two indicative nonlinearities are shown in Figure~\ref{Fig: comparison of nonlinearities} and {their effect to the closed-loop system} compared to {that} of the previous linear controller. Observe that the nonlinear control curves \eqref{eq: nonlinear control} result in a {qualitatively comparable} closed-loop performance as the linear control strategies, {while additionally} enforcing injection constraints and reducing the control effort due to the deadzone.

{A formal comparison between the considered control laws in terms of specific closed-loop-performance metrics as in \cite{ET:16} requires additional theoretical analysis and structured numerical experiments that we leave for future work.}

\section{Conclusion and outlook} \label{sec:conclusion}

\textit{Summary}: We have proposed a novel frequency control approach that achieves both local asymptotic stability of the closed-loop equilibria of power systems, modeled as a nonlinear, differential-algebraic, dynamical system, and economic-dispatch optimality. {The control architecture is based on a semi-decentralized gather-and-broadcast protocol. Hence the communication requirements are significantly lower than those of distributed architectures, and we avoid certain shortcomings of the proposed decentralized or distributed architectures.}


\textit{Open problem}: Extensive numerical tests indicate that closed-loop local asymptotic stability holds true even if Assumption \ref{ass: coeff} on the cost functions is violated. Proving or disproving such a claim is currently an open problem. {Likewise, 
it is of interest to remove the Standing Assumption \ref{ass:convexity} to allow for non-differentiable or non-strictly convex cost functions that result in a wider range of admissible control strategies.}

\textit{Outlook}: An important extension would be the inclusion of forecasts and inter-temporal constraints into the semi-decentralized frequency control architecture, with the aim of designing predictive control actions, while maintaining minimal communication requirements. 
{ 
Note that market mechanisms inspired our gather-and-broadcast strategy, but its final implementation is a distributed feedback control based on frequency measurements. 
An interesting avenue for future research is to explicitly consider bidding schemes and market mechanisms that are interfaced with the continuous-time power system dynamics at periodic time instants.
}

\section*{Acknowledgements}
The authors want to thank Felix Kottmann for insightful discussions on market mechanisms and their distinctions and Aaron Lelouvier for extensive simulation studies. 


\balance

\bibliographystyle{unsrt}
\bibliography{./bib/library,./bib/freq_reg}

\end{document}